\documentclass[12pt,a4paper,fleqn]{article}

\usepackage{amsmath,amsfonts,amsthm,amssymb,color}
\usepackage{latexsym}
\usepackage[latin1]{inputenc}
\usepackage{amsmath}
\usepackage{amsfonts}
\usepackage{oldstyle}
\usepackage{stmaryrd}
\usepackage{color}
\usepackage{delarray}
\usepackage{eufrak}

\newtheorem{prop}{Proposition}[section]

\newtheorem{lemma}[prop]{Lemma}

\newtheorem{theorem}[prop]{Theorem}

\renewcommand{\geq}{\geqslant}
\def\leq{\leqslant}

\newcommand{\N}{\mathbb{N}}

\newcommand{\R}{\mathbb{R}}

\newcommand{\cl}{\mathcal L}

\newcommand{\lp}{\left(}
\newcommand{\rp}{\right)}
\newcommand{\lc}{\left[}
\newcommand{\rc}{\right]}
\newcommand{\lcl}{\left\{}
\newcommand{\rcl}{\right\}}
\newcommand{\lln}{\left|}
\newcommand{\rrn}{\right|}

\def\HH{\EuFrak H}

\def\1{{\mathbf{1}}}

\def\1{{\mathbf{1}}}
\def\0.5{{\frac{1}{2}}}

\def\H{\EuFrak H}

\begin{document}

\title{Multivariate central limit theorems for averages of fractional Volterra processes and applications to parameter estimation}

\author{Ivan Nourdin\footnote{Email: {\tt  ivan.nourdin@uni.lu}; IN was partially supported by the Grant F1R-MTH-PUL-15CONF. (CONFLUENT) from Luxembourg University}, \,David Nualart\footnote{Email: {\tt  nualart@ku.edu}; DN was partially supported by the NSF grant  DMS1208625.} \,  and Rola Zintout\footnote{Email: {\tt rola.zintout@univ-lorraine.fr}. }}
\date{}

\maketitle
\noindent
{\bf Abstract}.
The purpose of this paper is to establish the multivariate normal convergence for the average of certain
Volterra processes constructed from a  fractional Brownian motion with Hurst parameter $H>\frac12$.  Some applications to parameter estimation are then discussed.

\section{Introduction}
Let $B^H$ be a fractional Brownian motion with Hurst parameter $H>\frac12$.
In this paper, we deal with fractional Volterra processes $X_i$, $i=1,\ldots,k$, of the form
\begin{equation}\label{X}
X_i(t) =\int _{0}^t  x_i(t-s) dB^H(s),\quad t\geq 0,
\end{equation}
where $x_i:[0,\infty)\to\R$ are measurable functions satisfying suitable integrability conditions (to be precised later on).

The special case of $k=1$ and $x_1(u)=\sigma\,e^{-\theta u}$, with $\sigma,\theta>0$, corresponds to the fractional Ornstein-Uhlenbeck process, which may be defined as the unique solution to the stochastic differential equation
\begin{equation}\label{fOU}
\left\{
\begin{array}{ll}
\dot{X}(t)&=-\theta X(t) + \sigma \dot{B}^H(t),\quad t>0\\
X(0)&=0.
\end{array}\right.
\end{equation}
In (\ref{fOU}) and everywhere else, the dots over $X$ and $B^H$ are used to indicate differentiation with respect to $t$.
More generally, consider the following $p$-th order stochastic differential equation driven by $B^H$:
\begin{equation}\label{CARintro}
\left\{
\begin{array}{ll}
X^{(n)}(t)&=\sum_{j=0}^{k-1}\theta_{j} X^{(j)}(t)+\sigma \dot{B}^H(t), \quad t>0\\
X(0)&=\ldots=X^{(p-1)}(0)=0.
\end{array}\right.
\end{equation}
In (\ref{CARintro}), the superscript ${}^{(j)}$ denotes $j$-fold differentiation with respect to $t$.
Then, it can be proved that the $p$-dimensional process $(X, X^{(1)}, \dots, X^{(p-1)})$ is of the form (\ref{X}) with $k=p$  for suitable functions $x_i$, $i=1,\dots, k$.

Considering that data are sampled from an underlying continuous-time process,
the solution $X$ to (\ref{CARintro}) can, for instance, serve as a model for (possibly irregularly spaced) discrete time long memory data. In such a situation, parameters $\theta_0,\ldots,\theta_{k-1}$ are usually   unknown and,  therefore, they must be accurately calibrated from the observation of
$X$. This is how we naturally arrive to the issue of showing a central limit theorem (CLT) for the  parameter estimators in the model (\ref{CARintro}), in order, e.g., to construct confidence intervals.
For instance, in the situation where those estimators are obtained by means of the method of moments,   we may be naturally led   to show  a multivariate CLT for random vectors taking the form
\[
\left\{ \frac1{\sqrt{T}}\int_0^T \left[ \big(X^{(i_j)}(t)\big)^{m_j} - E\big(X^{(i_j)}(t)\big)^{m_j} \right]dt:\, j=1,\ldots,k \right\},
\]
for given powers $m_1,\ldots,m_k\in\N^*$ and differentiation indices $i_1,\ldots,i_k\in\{0,\ldots,p-1\}$.
We refer to \cite{Hu-Nualart,KL,MN}  for different types of   central limit theorems for parameter estimators in this type of models.

Motivated by these statistical  problems, our purpose is to derive general central limit theorems for functionals of the process $X$ solution to (\ref{X}). More precisely,  let $f_i:\R\to\R$, $i=1,\ldots,k$, be real mesurable functions satisfying
\begin{equation}\label{integr-intro}
\int_\R f_i(x)e^{-x^2/2}dx=0\quad\mbox{and}\quad \int_\R f_i^2(x)e^{-x^2/2}dx<\infty.
\end{equation}
The second condition in (\ref{integr-intro}) ensures that $f_i$ can be expanded in Hermite polynomials, namely
\begin{equation}\label{Hermite-f-intro}
f_i=\sum_{l=0}^\infty a_{i,l}H_l\quad \mbox{with $\sum_{l=0}^\infty l!a_{i,l}^2<\infty$},
\end{equation}
whereas from the first one we deduce that $a_{i,0}=0$.
The first goal of the present paper is to answer the following question.

\bigskip

{\bf Question A}:  As $T\to\infty$, can we exhibit reasonable conditions ensuring that a multivariate CLT holds for the family of random vectors $U_T=(U_{1,T},\ldots,U_{k,T})$? Here
\begin{equation}\label{U}
U_{i,T}=\frac{1}{\sqrt{T}}\int_0^T f_i\left(\frac{X_i(t)}{\sigma_i(t)}\right)dt,\quad \mbox{with $\sigma_i(t)=\sqrt{E[X_i(t)^2]}$,}
\end{equation}
and $X_i$. $i=1, \dots, k$ are the fractional Volterra processes solution to  (\ref{X}).

\bigskip

Since, in general,  the processes $X_i$ are   not stationary, we stress that one cannot directly apply the classical Breuer-Major theorem (see Theorem \ref{BM}) to positively answer Question A. Nevertheless, following the approach developed in Nourdin, Peccati and Podolskij   \cite{NPP} (see also \cite[Chapter 7]{NPBook}) we prove the following result.

\begin{theorem}\label{main1}
Let $q_i$ denote
the {\sl Hermite rank} of $f_i$, that is, the smallest value of $l$ such that the coefficient $a_{i,l}$ of $H_l$ in (\ref{Hermite-f-intro}) is different from zero.
Set $q_*=\min_{1\leq i\leq k}q_i$ and assume that $q_*\geq 2$.
Consider $U_T=(U_{1,T},\ldots,U_{k,T})$, where $U_{i,T}$ is given by (\ref{U}).
If $H\in (\frac12,1-\frac{1}{2q_*})$ and if the functions $x_i$ defining $X_i$ satisfy
both \begin{equation}\label{breuer-cond}
\int_{0}^\infty
\left(\int_{[0,\infty)^2}|x_i(u)x_j(v)|\,|v-u-a|^{2H-2}dudv\right)^{q_i\vee q_j} da<\infty
\end{equation}
and
\begin{equation}\label{positivity}
\eta_{i} := \sqrt{H(2H-1)\int_{[0,\infty)^2} x_i(u)x_i(v)|v-u|^{2H-2}dudv}\in (0,\infty),
\end{equation}
for all $i,j=1,\dots, k$, then
\begin{equation}\label{conclusion1}
U_T
\overset{\rm law}{\to} N_k(0,\Lambda)\quad\mbox{ as $T\to\infty$},
\end{equation}
where $\Lambda=(\Lambda_{ij})_{1\leq i,j\leq k}$ is given by
\begin{eqnarray}\label{Lambda-intro}
\Lambda_{ij}&=&\sum_{l=q_i\vee q_j}^{\infty} a_{i,l} a_{j,l} l !  \frac{H^l (2H-1)^l }{\eta_i^l\eta_j^l}\\
&&\hskip1cm\times  \int_{\R} \lp \int_{[0,\infty)^2} x_i(u)x_j(v)|v-u-a|^{2H-2}dudv\rp^l da.\notag
\end{eqnarray}
\end{theorem}

That one must divide by a quantity depending on $t$ in (\ref{U}), namely $\sigma_{i}(t)$,  may appear to be not very convenient for applications. This is why we also address the following related problem.

\bigskip

{\bf Question B}: Can one find  constants $\xi_{i}>0$, as well as suitable assumptions on $f_i$, $x_i$ and $H$, so that
 $V_T=(V_{1,T},\ldots,V_{k,T})$ satisfies a CLT? Here
\begin{equation}\label{V}
V_{i,T}=\frac{1}{\sqrt{T}}\int_0^T f_i\left(\frac{X_{i}(t)}{\xi_{i}}\right)dt.
\end{equation}

\bigskip

Whatever the value of $\xi_i$, observe that the variance of $X_{i}(t)/\xi_{i}$ is  different from 1 for most of the values of $t$. For this reason and because Hermite polynomials are the orthogonal polynomials associated with the {\it standard} Gaussian distribution, it seems difficult to deal with {\it general} functions $f_i$ while trying to answer Question B. This is why we restrict our analysis to the situation where $f_i$ are {\it polynomials}, which is not a loss of generality for the applications we have in mind (see Section \ref{sec-car}).
More precisely, we have the following result, which provides a positive answer to Question B.

\begin{theorem}\label{main2}
 Suppose that    $f_i=P_i$, $i=1,\ldots,k$, are real polynomials
 and  denote by $q_i$ the Hermite rank of $P_i$. Set $q_*=\min_{1\leq i\leq k}q_i$ and assume that $q_*\geq 2$.
Consider $V_T=(V_{1,T},\ldots,V_{k,T})$ given by (\ref{V}),
where $\xi_i=\eta_i$ is given by (\ref{positivity}).
If $H\in (\frac12,1-\frac{1}{2q_*})$ and if the functions $x_i$ defining $X_i$ satisfy
 (\ref{breuer-cond}), (\ref{positivity}) as well as
 \begin{equation}\label{dol2-intro}
\int_{[0,\infty)^2} |x_i(u)x_i(v)|\,\big((u\wedge v) \vee1\big)\,|v-u|^{2H-2}dudv<\infty,
\end{equation}
 then
\begin{equation}\label{conclusion3}
V_T
\overset{\rm law}{\to} N_k(0,\Lambda)\quad\mbox{ as $T\to\infty$},
\end{equation}
with $\Lambda$ still given by (\ref{Lambda-intro}).
\end{theorem}

In the last section of our paper, we discuss an application of Theorem \ref{main2} to the problem of parameter estimation in the fractional  CAR($k$) model, which generalizes the model introduced in \cite{MN}.\\

The paper is organized as follows. Section 2 contains preliminary results and concepts. The proof of the two main results, namely Theorems \ref{main1} and \ref{main2}, is then provided in Section 3. Finally, an application of Theorem \ref{main2} is discussed in Section 4.

\section{Preliminaries}

\subsection{Fractional Brownian motion}

Throughout the paper, $B^H=(B^H(t))_{t\in\R}$ denotes a fractional Brownian motion (fBm in short) with  Hurst index $H\in (\frac12,1)$, defined on a complete probability space $(\Omega, \mathcal{F}, P)$. That is, $B^H$ is a zero mean Gaussian process with covariance
\[
E[B^H_t B^H_s] = \frac 12 \left( |t|^{2H} + |s| ^{2H} - |t-s|^{2H} \right).
\]
We further assume that the $\sigma$-field $\mathcal{F}$ is the completion of the $\sigma$-field generated by $B^H$.
We denote by $\H$  the closure of the space of step functions on $\R$ endowed with the inner product
\[
\langle \mathbf{1}_{[a,b]} ,  \mathbf{1}_{[c,d]} \rangle_\H = E[(B^H_b-B^H_a)(B^H_d-B^H_c)],
\]
for any $a<b$ and $c<d$.  We know that the Hilbert space $\H$ is isometric to the Gaussian space spanned by $B^H$, and we denote this isometry equivalently by $x\rightarrow B^H(x)$ or by $x\mapsto \int_{-\infty}^\infty x(s) dB^H(s)$.
\subsection{Wiener integral against fBm}
It is well-known (see, for instance, Nualart \cite[Chapter 5]{Nualart}) that any measurable function $x:\R\rightarrow \R$ satisfying
\begin{equation}\label{condition}
\int_{\R^2} |x(u)x(v)||v-u|^{2H-2}dudv <\infty
\end{equation}
belongs to the space $\H$, that is, it can be integrated with respect to $B^H$.
In this case,
the Wiener integral $B^H(x)=  \int_{-\infty}^\infty x(s) dB^H(s)$ satisfies the following isometry property:
\[
E \left[\left|    B^H(x) \right|^2\right]=  \|x \|^2_\H = H(2H-1) \int_{\R^2} x(u)x(v) |v-u|^{2H-2}dudv.
\]
Notice that condition (\ref{condition}) can be equivalently rewritten as
\begin{equation}\label{convol}
\int _\R ( |x|* |\widetilde{x}|)(t) |t|^{2H-2}  dt <\infty,
\end{equation}
with $\widetilde{x}(t)= x(-t)$ and where $x*y$ denotes the convolution of two nonnegative or integrable functions $x,y :\R\rightarrow \R$:
\[
(x*y) (t) =\int _{\R} x(u-v) y(v) dv.
\]

The following estimate will be needed in the proof of Theorem \ref{main1}.
\begin{lemma}\label{lemmaA}
For every $x,y:[0,\infty)\to\R_+$ satisfying the condition (\ref{condition})
(extended with $x(u)=y(u)=0$ whenever $u<0$) and for any $a\in\R$, we have
\begin{eqnarray*}
&&\left(\int_{[0,\infty)^2}x(u)y(v)|v-u-a|^{2H-2}dudv\right)^2\\
&\leq&
\int_{[0,\infty)^2}x(u)x(v)|v-u|^{2H-2}dudv
\int_{[0,\infty)^2}y(u)y(v)|v-u|^{2H-2}dudv.
\end{eqnarray*}
\end{lemma}
\noindent
{\it Proof}.
 We have, by Cauchy-Schwarz inequality,
\begin{eqnarray*}
&&\left(\int_{0}^\infty du\,x(u)\int_{0}^{\infty}dv\,y(v)|v-u-a|^{2H-2}\right)^2\\
&=&\left(\int_{0}^{\infty} du\,x( u)\int_{-a}^{\infty}dv\,y( v+a)|v-u|^{2H-2}\right)^2\\
&=&\frac{1}{H^2(2H-1)^2}\left(E\left[\int_{0}^{\infty} x( u)dB^H(u)\int_{-a}^{\infty}y( v+a)dB^H(v)\right]\right)^2\\
&\leq&\frac{1}{H^2(2H-1)^2}E\left[\left(\int_{0}^{\infty} x( u)dB^H(u)\right)^2\right]
E\left[\left(\int_{-a}^{\infty}y( v+a)dB^H(v)\right)^2\right]\\
&=&
\int_{[0,\infty)^2}x(u)x(v)|v-u|^{2H-2}dudv\int_{[0,\infty)^2}y(u)y(v)|v-u|^{2H-2}dudv.
\end{eqnarray*}
\qed

\subsection{Hermite polynomials and Wiener chaoses}
 For any integer $p\geq 1$, we denote by $%
\EuFrak H^{\otimes p}$ and $\EuFrak H^{\odot p}$, respectively, the $p$th tensor product and the $p$th symmetric tensor product of $\EuFrak H$.
   The  $p$th {\it Wiener chaos} of $X$, denoted by $\mathcal{H}_{p}$,  is  the closed linear subspace of $L^{2}(\Omega)$
generated by the random variables $\{H_{p}(B^H(x)),\,x\in \EuFrak H,\left\|
x\right\| _{\EuFrak H}=1\}$, where $H_{p}$ is the $p$th Hermite polynomial defined by
$$
H_p(x)=(-1)^p e^{x^2/2}\frac{d^p}{dx^p}\big(e^{-x^2/2}\big).
$$
The mapping $I_{p}(x^{\otimes p})=H_{p}(B^H(x))$ provides a
linear isometry between $\EuFrak H^{\odot p}$
(equipped with the modified norm $\sqrt{p!}\left\| \cdot \right\| _{\EuFrak %
H^{\otimes p}}$) and $\mathcal{H}_{p}$ (equipped with the $L^2(\Omega)$ norm).

\subsection{Fourth moment theorem}

The following result, known as the {\it fourth moment theorem}, is a combination of the seminal results of Nualart and Peccati \cite{NuPec} and Peccati and Tudor \cite{PecTud}. Given   a sequence of random vectors in a fixed Wiener chaos whose covariance matrix converges, the fourth moment theorem provides necessary and sufficient conditions for the convergence to a normal distribution.
We refer to Nourdin and Peccati  \cite{NPBook} for an extensive discussion  on this theorem, including  quantitative versions obtained by means of   Stein's method and a wide range of  applications and developments.

In the fourth moment theorem, a crucial role is played by the notion of {\it contractions}.
Let $\{e_{k},\,k\geq 1\}$ be a complete orthonormal system in $\EuFrak H$.
Given $f\in \EuFrak H^{\odot p}$, $g\in \EuFrak H^{\odot q}$ and $%
r\in\{0,\ldots ,p\wedge q\}$, the $r$th {\it contraction} of $f$ and $g$
is the element of $\EuFrak H^{\otimes (p+q-2r)}$ defined by
\begin{equation}
f\otimes _{r}g=\sum_{i_{1},\ldots ,i_{r}=1}^{\infty }\langle
f,e_{i_{1}}\otimes \ldots \otimes e_{i_{r}}\rangle _{\EuFrak H^{\otimes
r}}\otimes \langle g,e_{i_{1}}\otimes \ldots \otimes e_{i_{r}}\rangle _{%
\EuFrak H^{\otimes r}}.  \label{v2}
\end{equation}

\begin{theorem} [Fourth Moment Theorem] \label{4th}
Let $k \geq 2$ and $q_k\geq\ldots\geq q_1\geq 1$ be some fixed integers, and consider a family of kernels
\begin{eqnarray*}
\lcl  (f_{1,T},\ldots,f_{k,T})\rcl_{T>0}
\end{eqnarray*}
such that $f_{j,T}\in \HH^{\odot q_j}$ for every $T>0$ and every $j=1,\ldots,k$. Assume further that
\begin{eqnarray*}
\lim_{T \rightarrow\infty}E\lc I_{q_i}(f_{i,T}) I_{q_j}(f_{j,T})\rc&=&\Lambda_{ij},\quad \forall\mbox{ } 1\leq i,j\leq k.
\end{eqnarray*}
Then the following two conditions are equivalent:\\
(i) For every $i=1,\ldots,k$ and every $p=1,\ldots,q_i-1$,
\begin{eqnarray*}
\lim_{T \rightarrow\infty}\|f_{i,T}\otimes_p f_{i,T}\|_{\HH^{\otimes 2 (q_i-p)}}=0;
\end{eqnarray*}
(ii) as $T \rightarrow\infty$, the vector $\lp I_{q_1}(f_{1,T}),...,I_{q_k}(f_{k,T})\rp$ converges in distribution to the
$k$-dimentional Gaussian vector $N_k(0,\Lambda)$.
\end{theorem}


\subsection{Breuer-Major theorem}
We conclude this preliminary section with a continuous-time version of the celebrated Breuer-Major CLT for stationary Gaussian sequences.
\begin{theorem} [Breuer-Major]\label{BM}
Let $(X(t))_{t\geq  0}$ be a zero mean stationary Gaussian process with unit variance,
and let $f:\R\rightarrow \R$ be a measurable function
satisfying $\int_{\R}f^2(x)e^{-x^2/2}dx<\infty$.
Let us expand $f$ in terms of Hermite polynomials, namely
\[
f=\sum_{l=0}^\infty a_{l}H_l\quad \mbox{with $\sum_{l=0}^\infty l!a_{l}^2<\infty$}.
\]
Suppose that $a_0=0$ and let $q$ denote the Hermite rank of $f$ (that is, $q$ is the smallest value of $l$ such that the coefficient $a_{l}$ of $H_l$ is different from zero).
Finally, assume that  $\int_\R\lln \rho(t) \rrn^qdt<\infty$, with $\rho:\R\to\R$ the autocovariance function associated with $X$, that is,
\[
E[X(s)X(t)]=\rho(t-s),\quad t,s\geq 0.
\]
Then, as  $T\rightarrow\infty$,
\begin{eqnarray}
\frac{1}{\sqrt{T}}\int_0^T f(X(t))dt\mbox{ }\overset{\rm Law}{\to }N(0,\sigma^2),
\end{eqnarray}
where $\sigma^2=\sum_{l=q}^{\infty}l!a_l^2\int_{\R}\rho(t)^ldt \in(0,\infty)$.
\end{theorem}

\section{Proofs of the main results}

\subsection{Proof of Theorem \ref{main1}}

For the sake of clarity, the proof of Theorem \ref{main1} is divided into several steps.\\

\underline{Step 1}: {\sl Going away from zero}.
We claim that, in order to prove (\ref{conclusion1}) it is enough to show that,
for at least one {\it fixed} $T_0>0$, one has
\begin{equation}\label{conclusion2}
\widetilde{U}_T
\overset{\rm law}{\to} N_k(0,\Lambda)\quad\mbox{as $T\to\infty$},
\end{equation}
where $\widetilde{U}_T =(\widetilde{U}_{1,T},\ldots,\widetilde{U}_{k,T})$, with
\[
\widetilde{U}_{j,T} =
\frac{1}{\sqrt{T}}\int_{T_0}^T f_j\left(\frac{X_j(t)}{\sigma_j(t)}\right)dt
.
\]
(The only difference between $U_{j,T}$ and $\widetilde{U}_{j,T}$ is that the integral defining the former is between 0 and $T$.)
Indeed, using among other the Hermite expansion (\ref{Hermite-f-intro}) of each $f_j$ and then Cauchy-Schwarz inequality to obtain that $|E[X_{j}(t)X_{j}(s)]|\leq \sigma_j(t)\sigma_j(s)$, we get
\begin{eqnarray*}
&&E\left[\left(U_{j,T}-\widetilde{U}_{j,T} \right)^2\right]\\
&=&\frac1T \sum_{l=q_j}^\infty l!a_{j,l}^2\int_{[0,T_0]^2} \left( \frac{E[X_j(t)X_j(s)]}{\sigma_j(t)\sigma_j(s)}\right)^ldsdt\leq\frac{T_0^2}T \sum_{l=q_j}^\infty l!a_{j,l}^2=O(T^{-1}),
\end{eqnarray*}
as $T\to\infty$.
Hence, if (\ref{conclusion2}) holds true, then (\ref{conclusion1}) takes place as well.\\\\
\underline{Step 2}: {\sl Checking that the covariance matrix of $\widetilde{U}_T$ converges whenever $f_j$ are Hermite polynomials}.
 We can write, for any $j$ and as $t\to\infty$,
\begin{eqnarray*}
\sigma_j(t)^2 &=& H(2H-1)\int_{[0,t]^2}x_j(u)x_j(v)|v-u|^{2H-2}dudv\\
&\to& H(2H-1)\int_{[0,\infty)^2}x_j(u)x_j(v)|v-u|^{2H-2}dudv=\eta_j^2.
\end{eqnarray*}
Let us choose $T_0>0$ large enough so that $\sigma_j(t)^2\geq \frac12\eta_j^2$ for all $t\geq T_0$ and all $j\in\{1,\ldots,k\}$.
We shall check the convergence of the covariance matrix of $\widetilde{U}_T$
with this $T_0$ at hand and when $f_j=H_{p_j}$, $j=1,\ldots,k$, for given integers $p_j\geq q_j$. In this case, one has, for any $T>T_0$ and any $i,j\in\{1,\ldots,k\}$,
\begin{eqnarray}
&&E[\widetilde{U}_{i,T}\widetilde{U}_{j,T}]={\bf 1}_{\{p_i=p_j\}}\frac{p_i!}{T}\int_{[T_0,T]^2} \left(\frac{E[X_i(s)X_j(t)]}{\sigma_i(s)\sigma_j(t)}\right)^{p_i} dsdt\notag\\
&=&{\bf 1}_{\{p_i=p_j\}}\frac{p_i!}{T}\int_{T_0}^T db\int_{b-T}^{b-T_0} da \left(\frac{E[X_i(b)X_j(b-a)]}{\sigma_i(b)\sigma_j(b-a)}\right)^{p_i}\notag\\
&=&{\bf 1}_{\{p_i=p_j\}}\frac{p_i!}{T}\int_{T_0}^T db\int_{b-T}^{b-T_0} da\notag\\
&&\hskip1cm\times\left(\frac{H(2H-1)}{\sigma_i(b)\sigma_j(b-a)}\int_{[0,b]\times[0,b-a]}x_i(u)x_j(v)|v-u-a|^{2H-2}dudv\right)^{p_i}\notag\\
&=&{\bf 1}_{\{p_i=p_j\}}{p_i}!\int_{\R}\Phi(T,a)da,\label{bla}\,
\end{eqnarray}
where
\begin{eqnarray*}
&&\Phi(T,a)={\bf 1}_{\{|a|\leq T-T_0\}}\,\frac{1}T\int_{T_0\vee(a+T_0)}^{T\wedge (a+T)}\\
&&\hskip1.5cm \left(\frac{H(2H-1)}{\sigma_i(b)\sigma_j(b-a)}\int_{[0,b]\times[0,b-a]}x_i(u)x_j(v)|v-u-a|^{2H-2}dudv\right)^{p_i}db.
\end{eqnarray*}
Using first that $\sigma_i(t)^2\geq \frac12\eta_i^2$ for all $t\geq T_0$ and all $i$ (by definition of $T_0$) and then Lemma \ref{lemmaA},
we deduce that, for any $a\in\R$ and any $T>T_0$,
\begin{eqnarray}
 && {\bf 1}_{\{p_i=p_j\}}  |\Phi(T,a)|   \notag \\
 &\leq &  {\bf 1}_{\{p_i=p_j\}}    \left( \frac {2H(2H-1)} {\eta_i\eta_j}\int_{[0,\infty)^2}|x_i(u)x_j(v)|\,|v-u-a|^{2H-2}dudv\right)^{p_i}\notag\\
&\leq &  2^{p_i}     \left(  \frac{ H(2H-1) }{\eta_i\eta_j}\int_{[0,\infty)^2}|x_i(u)x_j(v)|\,|v-u-a|^{2H-2}dudv\right)^{q_i \vee q_j}.\notag\\
\label{blafar}
\end{eqnarray}
Moreover, due to the fact that $\sigma_i(t)\to\eta_i>0$, one has, as $T\to\infty$
\begin{eqnarray*}
&&\Phi(T,a)
={\bf 1}_{\{p_i=p_j\}}H^{p_i}(2H-1)^{p_i}{\bf 1}_{\{|a|\leq T-T_0\}}\int_{\frac{T_0}{T}\vee\frac{a+T_0}{T}}^{1\wedge (\frac{a}{T}+1)} \\
&&\hskip.5cm\left(\frac{\int_{[0,bT]\times[0,bT-a]}x_i(u)x_j(v)|v-u-a|^{2H-2}dudv}{\sigma_i(bT)\sigma_j(bT-a)}\right)^{p_i}db\\
&\to&{\bf 1}_{\{p_i=p_j\}} \left(
\frac{H(2H-1)}{\eta_i\eta_j} \int_{[0,\infty)^2}x_i(u)x_j(v)|v-u-a|^{2H-2}dudv\right)^{p_i}.
\end{eqnarray*}
By dominated convergence (see also (\ref{breuer-cond})), one deduces that
\begin{eqnarray}
E[\widetilde{U}_{i,T}\widetilde{U}_{j,T}]&\to& {\bf 1}_{\{p_i=p_j\}} \frac{p_i!H^{p_i}(2H-1)^{p_i}}{\eta_i^{p_i}\eta_j^{p_i}}\label{pc}\\
&&\times \int_{\R}  \left(\int_{[0,\infty)^2}x_i(u)x_j(v)|v-u-a|^{2H-2}dudv\right)^{p_i} da.\notag
\end{eqnarray}
Observe that (\ref{pc}) coincides with $\Lambda_{ij}$ after taking into account that $f_i=H_{p_i}$ and $f_j=H_{p_j}$.\\

\underline{Step 3}: {\sl Proving (\ref{conclusion2}) whenever $f_j$ are Hermite polynomials}.
To do so, we shall make use of the Fourth Moment Theorem \ref{4th}.
As in the previous step, suppose that $f_j=H_{p_j}$, $j=1, \dots ,k$, for some $p_j \geqslant q_j$. One then has
$\widetilde{U}_{j,T}=I_{p_j}\lp g_{j,T}\rp$, with
\begin{equation}\label{gjt}
g_{j,T}=\frac{1}{\sqrt{T}}\int_{T_0}^{T}\frac{e_{j,t}^{\otimes p_j}}{\sigma_j(t)^{p_j}} dt.
\end{equation}
In (\ref{gjt}), $e_{j,t}$ is a short-hand notation for the function $u\mapsto x_j(t-u){\bf 1}_{[0,t]}(u).$\\
According to Theorem \ref{4th}, to conclude that (\ref{conclusion2}) takes place
we are left to check that, for all $j\in\{1,\dots ,k\}$ and all $r\in\{1, \dots ,p_j-1\}$,
\begin{eqnarray}\label{toprove}
\|g_{j,T}\otimes_r g_{j,T}\|\to0 \mbox{ as } T\to\infty.
\end{eqnarray}
Let us compute $g_{j,T}\otimes_r g_{j,T}$. We find
\begin{eqnarray*}
g_{j,T}\otimes_r g_{j,T}=\frac{1}{T}\int_{[T_0,T]^2} \frac{e_{j,t_1}^{\otimes (p_j-r)}\otimes e_{j,t_2}^{\otimes (p_j-r)}}{\sigma_j(t_1)^{p_j}\sigma_j(t_2)^{p_j}}E[X_j(t_1)X_j(t_2)]^r dt_1dt_2.
\end{eqnarray*}
As a result, using moreover that $\sigma_j(t)^2\geq \frac12\eta_j^2$ for $T>T_0$ and introducing
$Z_j(t)=\int_{-\infty}^t |x_j(t-s)|dB^H_s$ (after extending $B^H$ to the whole $\R$), we obtain
\begin{eqnarray}
&&\|g_{j,T}\otimes_r g_{j,T}\|^2\notag\\
&=&\frac{1}{T^2}\int_{[T_0,T]^4}\frac{1}{(\sigma_j(t_1)\sigma_j(t_2)\sigma_j(t_3)\sigma_j(t_4))^{p_j}} E[X_j(t_1)X_j(t_2)]^r E[X_j(t_3)X_j(t_4)]^r\notag\\
&&\hskip3cm \times E[X_j(t_1)X_j(t_3)]^{p_j-r} E[X_j(t_2)X_j(t_4)]^{p_j-r}dt_1dt_2dt_3dt_4\notag\\
&\leq&\frac{4^{p_j}}{\eta_j^{4p_j}\,T^2}\int_{[0,T]^4} E[Z_j(t_1)Z_j(t_2)]^r E[Z_j(t_3)Z(t_4)]^r\notag\\
&&\hskip3cm \times E[Z_j(t_1)Z_j(t_3)]^{p_j-r} E[Z_j(t_2)Z_j(t_4)]^{p_j-r}dt_1dt_2dt_3dt_4.\notag\\
\label{coucou}
\end{eqnarray}
Set $\tau_j^2=E[Z_j(0)^2]$ and let $h_{j,T}$ denote the function
\[
h_{j,T}=\frac{1}{\tau_j^{p_j}\sqrt{T}}\int_{0}^T \widetilde{e}_{j,t}^{\otimes p_j}dt,
\quad\mbox{where }\widetilde{e}_{j,t}(u)=|x_j(t-u)|{\bf 1}_{(-\infty,t]}(u).
\]
Expressing the right-hand side of (\ref{coucou}) by means of $\| h_{j,T}\otimes_r h_{j,T}\|^2$ leads to
\begin{equation}\label{in}
\|g_{j,T}\otimes_r g_{j,T}\|^2\leq
\frac{4^{p_j}\tau_j^{4p_j}}{\eta_j^{4p_j}}
\|h_{j,T}\otimes_r h_{j,T}\|^2.
\end{equation}
On the other hand, it is straightforward to check that $Z_j$ is a {\it stationary} Gaussian process
and that
\[
I_{p_j}(h_{j,T})=\frac{1}{\sqrt{T}}\int_0^T H_{p_j}\left(\frac{Z_j(t)}{\tau_j}\right)dt.
\]
From assumption (\ref{breuer-cond}) and with Lemma 2.1, we deduce that $\int_{\mathbb{R}} | \rho_{Z_j}(t) | ^{p_j}  dt<\infty$. Indeed,
\[
\rho_{Z_j}(t)= H(2H-1) \int_{[0,\infty)^2} |x_j(u) x_j (v)| |u-v-t|^{2H-2} dudv.
\]
As a consequence, Breuer-Major Theorem \ref{BM} implies that
$I_{p_j}(h_{j,T})$ converges in law to a Gaussian. According to Theorem \ref{4th}, one
deduces that $\|h_{j,T}\otimes_r h_{j,T}\|\to 0$ as $T\to\infty$, implying in turn
that (\ref{toprove}) holds true (see (\ref{in})), and thus completing the proof of (\ref{conclusion1}) in the particular case where $f_j=H_{p_j}$.\\

\underline{Step 4}: {\sl Proving (\ref{conclusion1}) whenever $f_j$ are polynomials}. More precisely, let us suppose in this step that, for each $j=1,\dots, k$ one has
\begin{equation}\label{fj}
f_j = \sum_{l=q_j}^{m_j} a_{j,l}H_l,
\end{equation}
for a {\it finite} integer $m_j$.
Owing to the Cram\'er-Wold device, it is actually immediate to apply the conclusion of Step 3 in order to get (\ref{conclusion1}) in the case where $f_j$ given by (\ref{fj}).\\

\underline{Step 5}: {\sl Proving (\ref{conclusion1}) in all generality}.
Finally, let us consider the general situation of $f_j$ given by (\ref{Hermite-f-intro}).
To reach the conclusion in this case, and taken into account Step 4, it remains to show that, for any fixed $j=1,\ldots,k$,
\[
\frac{1}{\sqrt{T}}
\int_{T_0}^T
\sum_{l=m}^\infty a_{j,l}H_l\left(
\frac{X_j(t)}{\sigma_j(t)}
\right)dt\to 0\quad\mbox{in $L^2$ as $m\to\infty$.}
\]
Using identity (\ref{bla}) and its associated bound (\ref{blafar}),
one has, for any fixed $j=1,\dots ,k$:
\begin{eqnarray*}
&&E\left[\left(\frac{1}{\sqrt{T}}
\int_{T_0}^T
\sum_{l=m}^\infty a_{j,l}H_l\left(
\frac{X_j(t)}{\sigma_j(t)}
\right)dt
\right)^2\right]\\
&=&\frac1T\int_{[T_0,T]^2}\sum_{l=m}^\infty l!a_{j,l}^2\left(\frac{E[X_j(t)X_j(s)]}{\sigma_j(t)\sigma_j(s)}\right)^ldsdt\\
&\leq&\frac1T\int_{[T_0,T]^2}\left(\frac{E[X_j(t)X_j(s)]}{\sigma_j(t)\sigma_j(s)}\right)^{q_j}dsdt\times \sum_{l=m}^\infty l!a_{j,l}^2\\
&\leq&\sum_{l=m}^\infty l!a_{j,l}^2\\
&& \times
\int_{\R} \left(\frac{2H(2H-1)}{\eta_j^{2}}\int_{[0,\infty)^2}|x_j(u)x_j(v)|\,|v-u-a|^{2H-2}dudv\right)^{q_j} da,
\end{eqnarray*}
which, thanks to (\ref{breuer-cond}) and Lemma 2.1, tends to zero as $m\to\infty$.
\qed

\subsection{Proof of Theorem \ref{main2}}

Set
\begin{eqnarray*}
Y_{j,t}=\int_{-\infty}^t x_j(t-s)dB^H(s)\quad\mbox{and}\quad L_{j,t}=\int_{-\infty}^0 x_j(t-s)dB^H(s),
\end{eqnarray*}
so that $Y_j=X_j+L_j$. It is straightforward to check that $Y_j$ is stationary and that $\eta_j^2={\rm Var}(Y_j(0))$.

We have $\eta_j^2={\rm Var}(Y_j(t))$ for all $t$ by stationarity. So, since $P$ has Hermite rank $q_j$ and since (\ref{breuer-cond}) takes place, by applying the same arguments than the ones used in the proof of Theorem \ref{main1} but with $Y_j$ instead of $X_j$, we can prove that, as $T\to\infty$,
\begin{eqnarray*}
\lp \frac{1}{\sqrt{T}} \int_0^T P_1\lp \frac{Y_1(t)}{\eta_1}\rp dt,\ldots,\frac{1}{\sqrt{T}} \int_0^T P_k\lp \frac{Y_k(t)}{\eta_k}\rp dt\rp
\overset{\rm law}{\to} N_k(0,\Lambda).
\end{eqnarray*}
Thus, to reach the desired conclusion it suffices to show that, under (\ref{dol2-intro}) and for any fixed $j=1,\ldots,k$, and for any integer $p\geq 1$, 
\begin{equation}\label{*}
\frac{1}{\sqrt{T}}
\int_0^T  (X_{j,t}^p - Y_{j,t}^p)dt\overset{L^1}\to 0\quad\mbox{as $T\to\infty$}.
\end{equation}
We shall divide the proof of (\ref{*}) into two steps.\\

\underline{Step 1}: We shall show that $E\int_0^\infty L_{j,t}^{2k} dt$ is bounded for any integer $k\geq 1$.
Indeed, one can write, with $\mu_{2k}$ denoting the $k$-th even moment of the standard Gaussian,
\begin{eqnarray*}
&&E\int_0^\infty L_{j,t}^{2k} dt \\
&=& \mu_{2k} H^{k}(2H-1)^{k}\int_0^\infty dt
\left(\int_{(-\infty,0]^2}dudv\, x_j(t-u)x_j(t-v)|v-u|^{2H-2}
\right)^k\\
&\leq&\mu_{2k} H^{k}(2H-1)^{k}\int_0^\infty dt
\left(\int_{[t,\infty)^2}dudv\, | x_j(u)x_j(v)||v-u|^{2H-2}
\right)^k\\
&\leq&\mu_{2k} H^{k}(2H-1)^{k}
\int_{[0,\infty)^{2k}}du_1\ldots du_{2k}\, |x_j(u_1)|\ldots|x_j(u_{2k})|\\
&&\times|u_2-u_1|^{2H-2}
\ldots|u_{2k}-u_{2k-1}|^{2H-2}
\times \min\{u_i,\,i=1,\ldots,2k\}.
\end{eqnarray*}
Now, using that $\min\{u_i,\,i=1,\ldots,2k\}\leq \sum_{i=1}^k \min\{u_{2i-1},u_{2i}\}$, and taking into account condition
(\ref{dol2-intro}), we deduce
\begin{eqnarray*}
&&E\int_0^\infty L_{j,t}^{2k} dt \\
&\leq& k\mu_{2k} H^{k}(2H-1)^{k}
\left(\int_{[0,\infty)^{2}}dudv\, |x_j(u)x_j(v)||v-u|^{2H-2}\right)^{k-1}\\
&&\hskip3.5cm\times
\int_{[0,\infty)^{2}}dudv\,u\wedge v\, |x_j(u)x_j(v)||v-u|^{2H-2}\\
&\leq& k\mu_{2k} H^{k}(2H-1)^{k}
\left(\int_{[0,\infty)^{2}}dudv\, \big((u\wedge v)\vee 1\big)|x_j(u)x_j(v)||v-u|^{2H-2}\right)^{k}<\infty.\\
\end{eqnarray*}

\underline{Step 2}. Let us observe that
\begin{equation}
\left|\frac{1}{\sqrt{T}}
\int_0^T  (X_{j,t}^p - Y_{j,t}^p)dt\right|\leq \sum_{k=1}^p\binom{p}{k}\frac{1}{\sqrt{T}}\int_0^T dt|Y_{j,t}|^{p-k}|L_{j,t}|^k.
\end{equation}
For any fixed $p\geq 1$ and any $1\leq k\leq p$, one has, using the Cauchy-Schwarz inequality, Step 1 and that $Y$ is stationary,
\begin{eqnarray*}
&&\frac{1}{\sqrt{T}}E\int_0^T |Y_{j,t}|^{p-k}|L_{j,t}|^kdt\\
&=&\frac{1}{\sqrt{T}}E\int_0^{\rho T} |Y_{j,t}|^{p-k}|L_{j,t}|^kdt+\frac{1}{\sqrt{T}}E\int_{\rho T}^T |Y_{j,t}|^{p-k}|L_{j,t}|^kdt\\
&\leq& {\rm cst}\left(\sqrt{\rho}+\sqrt{\int_{\rho T}^\infty E[L_{j,t}^{2k}]dt} \right),
\end{eqnarray*}
where $0<\rho <1$. So, by letting $T\to\infty$ and then $\rho\to 0$, the desired conclusion (\ref{*}) follows, thus concluding the proof of Theorem \ref{main2}.\qed

\bigskip

\section{An application to the estimation of parameters in the fractional ${\rm CAR}(k)$ model}\label{sec-car}

Consider the fractional ${\rm CAR}(k)$ model, that is, the solution $X$ to:
\begin{equation}\label{CAR}
X^{(k)}(t)=\sum_{i=0}^{k-1}\theta_{i} X^{(i)}(t)+\sigma \dot{B}^H(t), \quad t>0.
\end{equation}
Here, $X^{(k)}$ indicates the $k$th derivative of the solution process $X$
and $\theta_i$ are real parameters considered as being unknown.
Moreover, up to appropriate scaling it is not a loss of generality to assume that $\sigma=1$.

In the sequel, we are going to illustrate  the use of our Theorem \ref{main2}
to the estimation problem in (\ref{CAR}). To keep the things as simple as possible, we shall only consider the case $k=2$, we shall assume zero initial conditions for $X$ and we shall put some restrictions on $\theta_0,\theta_1$.

More precisely, let $X$ be defined as the unique solution
to
\begin{equation}\label{CAR2}
\ddot{X}(t)=\theta_0 X(t)+\theta_1 \dot{X}(t) + \dot{B}^H(t),\quad X(0)=\dot{X}(0)=0,
\end{equation}
with $\theta_0,\theta_1<0$ and $\theta_1^2+4\theta_0 >0$.
The roots of the characteristic equation $r^2-\theta_1 r-\theta_0=0$ are
\begin{equation}\label{eq4.2}
p=\frac{\theta_1+\sqrt{\theta_1^2+4\theta_0}}{2} \quad\mbox{and}\quad q=\frac{\theta_1-\sqrt{\theta_1^2+4\theta_0}}{2}.
\end{equation}
The solution $X$ to (\ref{CAR2}) is given
by
\[
X(t)=\int_0^t  \frac{e^{p(t-s)}-e^{q(t-s)}}{p-q}dB^H(s) .
\]
The processes  $X_1=X$ and $X_2=\dot{X}$ are of the form (\ref{X}), with the corresponding functions
\[
x_1(t)= \frac { e^{pt} -e^{qt}}{p-q}
\quad\mbox{and}\quad
x_2(t)=\dot{x}_1(t)= \frac { pe^{pt} -qe^{qt}}{p-q}.
\]
We shall apply Theorem \ref{main2} to $X_1=X$ and $X_2=\dot{X}$ and to the polynomials
$P_1(x)= P_2(x) =x^2-1$, which have Hermite rank $1$.   
 Since $p$ and $q$ are negative numbers, the  functions $x_!$ and $x_2$  satisfy conditions (\ref{breuer-cond}) (\ref{positivity}) and (\ref{dol2-intro}).
As a consequence, Theorem \ref{main2} implies the following convergence in law as $T$ tends to infinity, provided $H\in (\frac 12, \frac 34)$:
\begin{equation}\label{cvloiY}
\sqrt{T}\lp  \frac 1T \int_0^T (X(t)^2, \dot{X}(t)^2) dt -m_\infty \rp\overset{\rm \cl}{\to }N_2(0,\Lambda),
\end{equation}
where $m_\infty= (\eta_1^2, \eta_2^2)$, with $\eta_i$, $i=1,2$, defined in (\ref{positivity}), and $\Lambda$ is the covariance matrix appearing in (\ref{Lambda-intro}).
In order to explicitly compute $m_\infty$ and $\Lambda$, we shall use a Fourier transform approach.

\bigskip

\noindent
{\it Computation of $m_\infty$:}
The first component of the vector $m_\infty$ is given by
\[
\eta_1^2=H(2H-1) \int_{\R}  (x_1*\tilde{x}_1)(t)  |t|^{2H-2} dt.
\]
The Fourier transform of  $x_1$ is given by
\[
\mathcal{F}x_1(\xi)=  \frac 1{p-q}  \left( \frac 1{p-i\xi} -\frac  1{q-i\xi} \right).
\]
On the other hand, the Fourier transform of $|t|^{2H-2}$ is $\kappa_{2H-2} |\xi| ^{1-2H}$, for some constant $\kappa_{2H-2}$.
Therefore, using Plancherel theorem we can write
\begin{eqnarray*}
\eta_1^2 &=&  \frac {d_H}{(p-q)^2} \int_\R   \left|  \frac 1{p-i\xi} -\frac  1{q-i\xi} \right|^2 |\xi|^{1-2H} d\xi \\
&=&   \frac {d_H}{(p-q)^2} \int_\R   \left(  \frac 1{p^2+ \xi^2} +\frac  1{q^2+\xi^2}  -\frac 2{ pq +\xi^2} \right) |\xi|^{1-2H} d\xi  \\
&=& \frac {e_H} {(p-q)^2} \left( |p|^{-2H} + |q|^{-2H} - 2 (pq)^{-H} \right)
=  \frac {e_H( |p|^{-H} - |q|^{-H})^2}  {(p-q)^2},
\end{eqnarray*}
where $e_H= 2d_H \int_0^\infty \frac { \xi^{1-2H} d\xi}{ 1+\xi^2}$ and $d_H= \kappa_{2H-2}H(2H-1)$.

For the second component of $m_{\infty} $ we can write
\[
\eta_2^2=E\left[ \left( \int^{\infty} _0 \frac{pe^{pt}-qe^{qt}}{p-q}dB^H(t) \right)^2 \right]
=H(2H-1) \int_{\R}   (x_2*\tilde{x}_2)(t)   |t|^{2H-2} dt.
\]
The Fourier transform of  $x_2$ is given by
\[
\mathcal{F}x_2(\xi)=  \frac 1{p-q}  \left( \frac p{p-i\xi} -\frac  q{q-i\xi} \right).
\]
Therefore,
\begin{eqnarray*}
\eta_2^2 &=&  \frac {d_H}{(p-q)^2} \int_\R   \left|  \frac p{p-i\xi} -\frac  q{q-i\xi} \right|^2 |\xi|^{1-2H} d\xi \\
&=&   \frac {d_H}{(p-q)^2} \int_\R   \left(  \frac {p^2}{p^2+ \xi^2} +\frac  {q^2}{q^2+\xi^2}  -\frac {2pq}{ pq +\xi^2} \right) |\xi|^{1-2H} d\xi  \\
&=& \frac {e_H} {(p-q)^2} \left( |p|^{2-2H} + |q|^{2-2H} - 2 (pq)^{1-H} \right)
=  \frac {e_H( |p|^{1-H} - |q|^{1-H})^2}  {(p-q)^2}.
\end{eqnarray*}

 \medskip
 \noindent
 {\it Computation of  the matrix $\Lambda$:}
From  (\ref{Lambda-intro}), taking into account that  $P_1=P_2= H_2$, we obtain
\[
\Lambda_{ij}=  2H^2(2H-1)^2   \int_{\R^3} (x_i  *\tilde{x}_j)(u) (x_i  *\tilde{x}_j)  (v)   |u+a|^{2H-2}  | v+a| ^{2H-2} dudvda.
\]
We know  that
\[
\int_{\mathbb{R}} |u+a|^{2H-2}  | v+a| ^{2H-2} da= k_H |u-v| ^{4H-3},
\]
for some constant  $k_H$. Therefore, using again the Fourier transform, we can write
 \begin{eqnarray*}
\Lambda_{ij} &=& k_H  2H^2(2H-1)^2
\int_{\R^2} (x_i  *\tilde{x}_j)(u) (x_i  *\tilde{x}_j)  (v) |u-v|^{4H-3} dudv\\
&=&  a_H   \int_{\R}   |\mathcal{F}x_i(\xi)|^2  |\mathcal{F}x_j(\xi)|^2 |\xi|^{2-4H} d\xi,
\end{eqnarray*}
where $a_H= k_H  2H^2(2H-1)^2  \kappa_{4H-3}$.
The previous computations  lead  to the following expression for the components of the matrix $\Lambda$.
\begin{eqnarray*}
\Lambda_{1,1} &=&  \frac 1{(p-q)^4} \Big( \alpha_H (|p|^{-1-4H} + |q| ^{-1-4H} +4 (pq)^{-\frac 12-2H} )\\
&&
+\beta_H  \Big( \frac 2 {p^2-q^2} (|q|^{1-4H} -|p|^{1-4H} )
-\frac 4{p^2-pq} ( (pq)^{\frac 12-2H} -|p|^{1-4H}  )\\
&&\quad -\frac 4{q^2-pq} ( (pq)^{\frac 12-2H} -|q|^{1-4H}) \Big) \Big),
\end{eqnarray*}
\begin{eqnarray*}
\Lambda_{1,2} &=&  \frac 1{(p-q)^4} \Big( \alpha_H (|p|^{1-4H} + |q| ^{1-4H} +4 (pq)^{\frac 12-2H} )\\
&&
+\beta_H  \Big( \frac {p^2+q^2} {p^2-q^2} (|q|^{1-4H} -|p|^{1-4H} )
-2\frac {p^2+pq}{p^2-pq} ( (pq)^{\frac 12-2H} -|p|^{1-4H}  )\\
&&\quad -2\frac {q^2+pq}{q^2-pq} ( (pq)^{\frac 12-2H} -|q|^{1-4H}) \Big) \Big),
\end{eqnarray*}
 and
\begin{eqnarray*}
\Lambda_{2,2} &=&  \frac 1{(p-q)^4} \Big( \alpha_H (|p|^{3-4H} + |q| ^{-4H} +4 (pq)^{-\frac 32-2H} )\\
&&
+\beta_H  \Big( \frac {2p^2q^2}{p^2-q^2} (|q|^{1-4H} -|p|^{1-4H} )
-\frac {4p^3q}{p^2-pq} ( (pq)^{\frac 12-2H} -|p|^{1-4H}  )\\
&&\quad -\frac {4pq^3}{q^2-pq} ( (pq)^{\frac 12-2H} -|q|^{1-4H}) \Big) \Big).
\end{eqnarray*}
In the above expressions the constants $\alpha_H$ and $\beta_H$ are given by
\[
\alpha_H= 2a_H \int_0^\infty \frac {\xi^{2-4H} d\xi}  {(1+ \xi^2)^2}
\quad
\mbox{and}
\quad
\beta_H= 2a_H \int_0^\infty \frac {\xi^{2-4H} d\xi}  {1+ \xi^2}.
\]


\end{document}